\documentclass[12pt]{article}
\usepackage{amsmath, amssymb,array,epsfig,amsthm}
\usepackage{graphicx}
\theoremstyle{plain}
\newtheorem{lemma}{Lemma}[section]

\newtheorem{theorem}[lemma]{Theorem}
\newtheorem{corollary}[lemma]{Corollary}

\theoremstyle{remark}
\newtheorem{remark}[lemma]{Remark}

\newcommand{\N}{{\mathbb N}}

\newcommand{\R}{{\mathbb R}}
\newcommand{\C}{{\mathbb C}}

\newcommand{\Dom}{{\rm Dom}}

\newcommand{\Spec}{{\rm Spec}}

\renewcommand{\Re}{\operatorname{Re}}
\renewcommand{\Im}{\operatorname{Im}}

\newcommand{\dsc}{\mathrm{disc}}
\newcommand{\ess}{\mathrm{ess}}

\newcommand{\point}{.}

\numberwithin{equation}{section}

\title{On Approximation of the Eigenvalues
of Perturbed Periodic Schr\"odinger Operators
\thanks{This research is supported by the Leverhulme Trust grant F/00 276/F.}}
\author{Lyonell Boulton and Michael Levitin\\
\normalsize\small Maxwell Institute for Mathematical Sciences\\
\normalsize\small and Department of Mathematics\\
\normalsize\small Heriot-Watt University\\
\normalsize\small Riccarton, Edinburgh EH14 4AS, U.~K.\\
\normalsize\small email {\sffamily \{L.Boulton, M.Levitin\}@ma.hw.ac.uk}\\
\normalsize\small {\sffamily www.ma.hw.ac.uk/\textasciitilde lyonell/}, {\sffamily www.ma.hw.ac.uk/\textasciitilde levitin/}}
\date{February 2007}

\begin{document}

\maketitle

\begin{abstract}
This paper addresses the problem of computing the eigenvalues
lying in the gaps of the essential spectrum of a periodic
Schr\"odinger operator perturbed by a fast decreasing potential.
We use a recently developed technique, the so called 
\emph{quadratic projection method},
in order to achieve convergence free from spectral pollution. 
We describe the theoretical foundations of the method in detail, 
and illustrate its effectiveness by several examples.
\end{abstract}

\section{Introduction}\label{sect:intro} %
It is well known that the problem of
approximating the  eigenvalues lying in gaps of the essential spectrum of a
self-adjoint operator by a sequence of finite-dimensional problems (e.g. 
for numerical analysis) is far from trivial. The presence of essential 
spectrum both above and below an eigenvalue 
means that there is no obvious  variational principle (cf. e.g.
\cite{des}), so an approximation/computation  by a standard projection
method is not always possible.  The main  difficulty is due to the existence of sequences of
eigenvalues  of the (finite-dimensional) approximate operators, accumulating at
points in the gaps which do not belong to the spectrum. 
These points are called \emph{spurious eigenvalues}, and the phenomenon itself is often referred to as
\emph{spectral pollution}.

It has been shown, for general unbounded self-adjoint operators,  that spectral
pollution in a projection method may occur at any real point of the resolvent 
set located between two
parts of the essential spectrum  (see \cite[Theorem~2.1]{lesh}). 
This is a consequence of the fact that the resolvent 
is not compact.  A substantial amount of research has been devoted to finding 
ways of
choosing the projectors,  in order to achieve a ``safe'' method for particular 
problems, 
see e.g. \cite{rasa} and \cite{bobr}.
Techniques vary considerably according to the problem and are by no means 
universal.

In this paper we address the question of spectral pollution and its avoidance 
for a perturbed periodic Schr\"odinger operator
\begin{equation} \label{eq:schr}
H:=-\Delta + V,
\end{equation}
acting in the Hilbert space $L^2(\R^N)$, where $V=V_p+V_d$, with $V_p$ being 
purely periodic with  respect to some
lattice of $\R^N$ and $V_d$ being fast decaying at infinity. The essential 
spectrum
of $H$ is determined by $V_p$. It consists of bands of absolutely continuous
spectrum, separated by gaps in the resolvent set. If $V_d=0$, the spectrum is
purely essential. If $V_d\not=0$, discrete eigenvalues may appear in the gaps, 
see \cite{dehe}.

A usual method for finding the essential spectrum of $H$
analytically, the so called Floquet-Bloch technique, has been 
well studied (see e.g. \cite{rs4}, \cite{ku} and the references therein). 
It gives a decomposition of the periodic part of the operator
into a direct integral
of operators on a basic periodic cell.
This reduces the problem of finding the endpoints of the bands in the
essential spectrum, to the problem of finding the eigenvalues of
differential operators in a compact domain with regular boundary conditions. 

Much less in known about the discrete spectrum of $H$, which has to be either
estimated numerically or studied by means of asymptotic techniques (for the
latter see  e.g. \cite{dehe}, \cite{bir} and \cite{sus}).  As we shall see
below, the natural approach of truncating $\R^N$ to a large compact domain  and
applying the projection method to the corresponding Dirichlet problem,  is 
prone
to spectral pollution. This makes the numerical localisation of
these eigenvalues particularly  difficult.

The purpose of this paper is to  describe an alternative procedure for finding 
eigenvalues, the so called 
\emph{quadratic projection method}, recently studied in an abstract setting
in \cite{shar}, \cite{lesh}, \cite{bo1} and \cite{bo2}.
The distinctive feature of our method is that the underlying
discretised eigenvalue problem is quadratic in the spectral parameter (rather than linear), and 
has non-real eigenvalues.
Its main advantage over a standard projection method lies in its robustness: 
it \emph{never} pollutes and it always provides 
\emph{a posteriori two-sided estimates} of the error of computed eigenvalues.

The  paper is organised as follows.
In Section~\ref{sect:quad} we discuss the phenomenon of spectral
pollution in a standard projection method and discuss the quadratic 
projection in an abstract context. Our Corollary~\ref{cor:est2} 
is an improvement upon previously known non-pollution results
for the general quadratic method. In Section~\ref{sect:schrod}
we provide details on how to implement the quadratic projection method 
for the numerical localisation of the eigenvalues of operator $H$. 
We also discuss some concrete numerical examples,
but deliberately avoid including  the full account of the numerical procedures we have used, 
in order not to overload the text with unnecessary technical
details. These will appear elsewhere.

\section{The quadratic projection method}\label{sect:quad}

\subsection{Spectral pollution in an ordinary projection method}\label{subsect:poll} 

Before proceeding to describe our method, we want to give a rigorous motivation
why it is needed at all, and why spectral pollution is intrinsic in the 
standard
projection method.

Let $A$  be a
self-adjoint operator in a Hilbert space
$\mathcal{H}$ with a  dense domain, $\mathrm{Dom}(A)$. The spectrum of  $A$, $\Spec(A)$,  can be 
decomposed into the discrete spectrum, $\Spec_\dsc(A)$, consisting of isolated eigenvalues of finite
multiplicity, and the essential spectrum, 
$\Spec_\ess (A):=\Spec(A)\setminus \Spec_\dsc(A)=\{\lambda:A-\lambda I\text{ is not Fredholm}\}$.

Take a finite-dimensional subspace $\mathcal{L}\subset  \mathrm{Dom}(A)$, and 
let $\Pi_\mathcal{L}:\mathcal{H}\longrightarrow \mathcal{L}$ be the orthogonal 
projection onto $\mathcal{L}$.
Let $A_\mathcal{L}:=\Pi_\mathcal{L}A\upharpoonright \mathcal{L}$.

The \emph{projection} method, also known as the \emph{Galerkin} method, 
consists in truncating the (infinite-dimensional) spectral problem
$A u=\lambda u$ to 
\begin{equation} \label{eq:proj}
   A_\mathcal{L}u=\lambda  u \qquad \mathrm{for\ some\ }
   u \in \mathcal{L}\setminus\{0\}.
\end{equation}    
If the operator $A$ is bounded from below and 
has a compact resolvent, this provides an effective way of estimating numerically 
the eigenvalues of $A$. 
The $k$-th eigenvalue of \eqref{eq:proj} will always be above
the  $k$-th eigenvalue of $A$, counting multiplicity, 
\cite[Section~XIII.1]{rs4}. 
Furthermore, if  $\mathcal{L}$ approximates $\Dom(A)$
reasonably well, then the first few 
eigenvalues of \eqref{eq:proj} will be close to the corresponding
ones of $A$.  

A precise statement can be easily obtained from the minimax principle: 
\begin{lemma} \label{lem:proj_conv}
Let $\mathcal{L}_n$ 
be a sequence of finite-dimensional subspaces of $\Dom(A)$.
Assume that $A$  is bounded below and has a compact resolvent.
Let $\lambda_1\leq \ldots
\leq \lambda_m$ be the first $m$ eigenvalues of $A$. Let 
\[
   \mathcal{E}=\mathrm{Span}\,\{ u\in \Dom(A): Au=\lambda_ku,\,
1\leq k \leq m \},
\]
be the spectral subspace associated with $\{\lambda_1,\dots,\lambda_m\}$. If
\begin{equation} \label{eq:conv_cond}
  \lim\limits_{n\to\infty} \|A^p(u-\Pi_{\mathcal{L}_n}u)\|= 0\,, 
\end{equation}   
holds for $p=0, 1$  and all $u\in \mathcal{E}$, then the $k$-th eigenvalue of \eqref{eq:proj} approaches
the $k$-th eigenvalue of $A$ as $n\to \infty$ for $1\leq k \leq m$. 
\end{lemma}
We omit the proof.

In some particular
cases it is also possible to estimate the convergence rate of the 
eigenvalues \cite{SF}.

Similar results can be established if the resolvent of $A$ is non-compact, 
for eigenvalues outside the extrema of $\Spec_\ess(A)$. However the situation 
changes if we want to approximate an eigenvalue in
a gap of $\Spec_\ess(A)$. There is no easy minimax principle, and spectral pollution may 
happen at any point of the gap. 

The difficulties involved in the computation of these eigenvalues  are 
well known for particular operators, see e.g.
\cite{bobr} or \cite{rasa}.  Moreover, in a generic situation we have
\begin{lemma} \label{lem:poll}
If $\lambda\not\in \Spec_\ess(A)$ is such that
$\alpha<\lambda<\beta$ where $\alpha,\beta\in \Spec_\ess(A)$, there exists
a sequence of subspaces $\mathcal{L}_n$ satisfying
\eqref{eq:conv_cond} for all $p\in\N$ and all $u\in \Dom(A)$, 
such that $\lambda\in \Spec(A_{\mathcal{L}_n})$ for all $n\in\N$.
\end{lemma}

This lemma directly follows from \cite[Theorem~2.1]{lesh}.

\subsection{The abstract quadratic projection method}\label{subsect:quad}

Let, as before, $\mathcal{L}$ be a finite-dimensional subspace of $\Dom(A)$, and
let $E=\{e_1,\dots,e_n\}$ be a basis of $\mathcal{L}$. This basis need not 
be orthogonal. 

Consider the quadratic matrix polynomial
\begin{equation} \label{eq:P}
     P_{\mathcal{L}}(z):=Q_{\mathcal{L}}-2 zA_{\mathcal{L}}+z^2 B_{\mathcal{L}},
\end{equation} 
where 
\begin{equation} \label{eq:BAQ}
    [B_{\mathcal{L}}]_{jk}=\langle e_j,e_k \rangle\,\qquad
    [A_{\mathcal{L}}]_{jk}=\langle Ae_j,e_k \rangle\,\qquad
    [Q_{\mathcal{L}}]_{jk}=\langle Ae_j,Ae_k \rangle\,.
\end{equation}
In numerical analysis, $A_{\mathcal{L}}$ is called the \emph{stiffness} matrix, 
$B_{\mathcal{L}}$ is a \emph{mass} matrix, and $Q_{\mathcal{L}}$ is a \emph{bending} matrix.
If $E$ is an orthonormal basis, then $A_{\mathcal{L}}=\Pi_\mathcal{L}A\upharpoonright \mathcal{L}$, and
$B_{\mathcal{L}}=\operatorname{Id}\upharpoonright \mathcal{L}$. Additionally, if $E\subseteq \Dom(A^2)$, then
$Q_{\mathcal{L}}=\Pi_\mathcal{L}A^2\upharpoonright \mathcal{L}$ and 
$P_{\mathcal{L}}(z)=\Pi_\mathcal{L}(A-z)^2\upharpoonright \mathcal{L}$.

We define the spectrum of the matrix polynomial $P_{\mathcal{L}}$, $\Spec(P_\mathcal{L})$, 
as the set of $\mu \in \C$ such that 
\begin{equation} \label{eq:specP}
   P_{\mathcal{L}}(\mu)u=0\qquad \text{for\ some } u\in \mathcal{L}\setminus\{0\}.
\end{equation} 
Since $B_{\mathcal{L}}$ is non-singular, $\det (P_{\mathcal{L}}(z))$ is a polynomial
in $z$ of degree $2\dim (\mathcal{L})$. Moreover, if $\mu\in \Spec (P_\mathcal{L})$,
then also $\overline{\mu}\in \Spec (P_\mathcal{L})$. Therefore 
$\Spec (P_\mathcal{L})$ is a set of at most $2\dim (\mathcal{L})$ complex points,
symmetric with respect to the real axis.

The core idea of the quadratic projection method lies in the fact that $\Spec(A)$ can be
well estimated if one knows the points of $\Spec(P_\mathcal{L})$ which are ``close'' to the real 
line, see Corollary~\ref{cor:est1}
and Theorem~\ref{thm:approx} below.
In \cite{shar}, \cite{lesh} and \cite{bo1},
$\Spec(P_{\mathcal{L}})$ is called the \emph{second order spectrum of $A$
relative to $\mathcal{L}$}. This set was first studied in
connection with the spectrum of $A$ in \cite{da0}, where the name
originated.

\begin{remark}
Intuitively, the quadratic projection method arises from the following simple observation. 
Let $\zeta\in \R$ lie in a gap of the essential spectrum. By virtue of
the spectral theorem, the discrete eigenvalues of $(\zeta-A)$
inside the corresponding shifted gap of $(\zeta-A)$ containing the origin, are also
the discrete eigenvalues of $(\zeta-A)^2$ lying below the bottom of the essential spectrum of
$(\zeta-A)^2$. This suggests that the truncations of the latter
operator must provide information about the localisation of a
portion of $\Spec_\dsc(A)$ near $\zeta$. 
The quadratic projection method  is a rigorous realisation of a similar idea.
\end{remark}

The main reason for preferring \eqref{eq:specP}
over \eqref{eq:proj} for estimating the spectrum of $A$ lies in the
following observation.
Let $D(a,b)$ be the open disk in the complex plane with an interval $[a,b]$ as a diameter:
\[
   D(a,b):=\left\{w\in \C\,:\, \left|w-\frac{a+b}{2}\right|<
   \frac{b-a}{2}\right\}.
\]

\begin{theorem}[{\cite[Lemma~5.2]{lesh}}]\label{thm:lesh}
Suppose that $(a,b)\cap \Spec(A)=\varnothing$. If $z\in D(a,b)$,
then the matrix $P_\mathcal{L}(z)$ is non-singular.
\end{theorem}
\begin{proof} Our proof is slightly different from that of \cite{lesh}. Let $z\not \in D(a,b)$. Let
\[
    \Sigma_z:=\{(\lambda-z)^2\,:\,\lambda\in(-\infty,a]\cup
    [b,\infty) \}.
\]
We first show that $0\not \in \overline{\mathrm{Conv}\, \Sigma_z}$
(here $\mathrm{Conv}\,\Omega$ denotes the convex hull of the set
$\Omega\subset \C$). Indeed, let $\theta$ be the angle at $z$ of
the triangle $T$ whose vertexes are $a,b,z$. Elementary
geometric arguments show that $\theta>\pi/2$. Then the
transformation $m:\lambda\mapsto (\lambda-z)^2$, maps the
angular region
\[
    B=\{(w-z)\,:\, \rho w\in T\ \mathrm{for\ some\ }\rho\geq 0\}
\]
into another angular sector centred at the origin with 
angle $2\theta>\pi$. Since $(-\infty,a)\cup(b,\infty)\subset
\C\setminus B$ and
\[
   m\,:\, (-\infty,a]\cup[b,\infty)\longmapsto \Sigma_z,
\]
there exists $-\pi< \theta_0\leq \pi$ and $c>0$, such that
$\mathrm{Re}\, (e^{i\theta_0} w)\geq c$ for all $w\in \Sigma_z$.
This ensures that $0\not \in \overline{\mathrm{Conv}\, \Sigma_z}$
as required.

Since $A=A^*$, $(A-z)^2$ with domain $\Dom (A^2)$
is a normal operator, \cite{kato}. As we have for the numerical range
\[
   \mathrm{Num}\,(A-z)^2 \subseteq \overline{
   \mathrm{Conv}\, [\Spec (A-z)^2]}
   \subseteq \overline{\mathrm{Conv}\, \Sigma_z}\,,
\]
and $\Dom (A^2)$ is a core for $A$, we have
\[
   \Re \big[ e^{i\theta_0}\langle (A-z)u,(A-\overline{z})u\rangle \big]
   \geq c
\]
for all $u\in \Dom (A)$ with $\|u\|=1$. In particular this
holds true for all $u\in \mathcal{L}$ with $\|u\|=1$, so that
$P_{\mathcal{L}}(z)$ cannot be a singular matrix. \end{proof}

As a consequence of Theorem~\ref{thm:lesh}, the points of $\Spec(P_\mathcal{L})$
which are close to the real line, are
\emph{necessarily} close to $\Spec(A)$. In other words, 
the method never pollutes. We also have two immediate corollaries.

\begin{corollary} \label{cor:est1}
If $\mu\in \Spec (P_\mathcal{L})$, then
\begin{equation} \label{eq:est1}
\inf\{|\Re \mu-\lambda|\,:\, \lambda\in \Spec (A)\}\leq |\Im \mu|.
\end{equation}
\end{corollary}

If $\lambda\in \Spec(A)$ is isolated 
from other point of the spectrum, 
\eqref{cor:est1} provides a two-sided estimate of $\lambda$, with an error explicitly
determined without the need for computing eigenfunctions.
In case this error is small, we can actually improve it by a square:

\begin{corollary} \label{cor:est2}
 Let $\lambda\in \Spec(A)$. Assume that $\lambda$ is
isolated from other points of the spectrum and let
\begin{equation} \label{eq:est2}
\begin{aligned}   
\delta&:=\min\{|\lambda-\nu|\,:\, \nu\not=\lambda,\, \nu \in \Spec (A)\}\\
     &=\mathrm{dist}\,(\lambda,\Spec(A)\setminus \{\lambda\}).
\end{aligned}
\end{equation}
If $|\mu-\lambda|<\delta/2$ for $\mu\in \Spec (P_\mathcal{L})$, then
\begin{equation} \label{eq:est3}
|\Re \mu -\lambda|< \frac{2(\Im \mu)^2}{\delta}.
\end{equation}
\end{corollary}
\begin{proof}
Theorem~\ref{thm:lesh} yields
\begin{equation} \label{eq:est4}
\left|\mu-\left(\lambda\pm\frac{\delta}{2}\right)\right|>\frac{\delta}{2}\,.
\end{equation}
Using the assumption $|\mu-\lambda|<\delta/2$, we can re-write \eqref{eq:est4} as 
\[
   |\Re \mu-\lambda|<\frac{\delta}{2}-\sqrt{\frac{\delta^2}{4}-(\Im \mu)^2}\,.
\]
Thus
\[
   |\Re \mu-\lambda|<\frac{(\Im \mu)^2}{\displaystyle\frac{\delta}{2}+\sqrt{\frac{\delta^2}{4}-(\Im \mu)^2}}
   <\frac{2(\Im \mu)^2}{\delta}\,.
\]
\end{proof}

Corollary~\ref{cor:est2} supersedes Corollary~\ref{cor:est1}
once we have found points of $\Spec (P_\mathcal{L})$ sufficiently close
to an isolated point of the spectrum of $A$. 
Note that $\lambda\in \Spec(A)$ does not have to be a discrete eigenvalue.

The above ``non-pollution'' results  are useful as long as there are points of
$\Spec(P_\mathcal{L})$ near to the real line. It is not immediately clear, however, whether or not
the eigenvalues of $A$ are approximated by some points in
$\Spec(P_\mathcal{L})$ when the dimension of $\mathcal{L}$ goes to infinity.
The results of \cite{bo1} and \cite{bo2} show that this is indeed the case, 
under a condition analogous to
\eqref{eq:conv_cond}.

\begin{theorem}[{\cite[Theorem~2.2]{bo2}}] \label{thm:approx} 
Let $\lambda \in \Spec_\dsc(A)$, and let 
$\tilde{\mathcal{E}}_\lambda:=\{u\,:\,Au=\lambda u\}$ be the corresponding eigenspace. 
Let $\mathcal{L}_n\subset \Dom(A^2)$ 
be subspaces with corresponding orthogonal projections
$\Pi_{\mathcal{L}_n}$, such that \eqref{eq:conv_cond} holds for
$p=0,1,2$ and all $u\in\tilde{\mathcal{E}}_\lambda$.
Then there exist eigenvalues $\lambda_n\in \Spec(P_{\mathcal{L}_n})$ such that $\lambda_n\to \lambda$
as $n\to\infty$.
\end{theorem}

\section{The quadratic projection method for 
perturbed periodic Schr\"odinger operators} \label{sect:schrod}


Let $H$ be the differential expression defined by \eqref{eq:schr} acting on
the dense domain $W^{2,2}(\R^N)$.

Let
\begin{gather*}
p=2\quad\text{ if }  N\leq 3, \\ 
p>2\quad\text{ if }  N=4,\\
p>N/2\quad\text{ if }  N\geq 5.
\end{gather*}
Below and elsewhere we assume that the potential
$V:\R^N\longrightarrow \R$ is \emph{uniformly locally $L^p$} in
the sense that
\begin{equation} \label{eq:V_cond}
   \int_{C} |V(x)|^p \,\mathrm{d} ^N x \leq M
\end{equation}
for any unit hyper-cube $C$, where the constant $M$ is independent
of $C$.

The condition \eqref{eq:V_cond}
ensures that the operator of multiplication by $V$
is $(-\Delta)$-bounded with relative bound equal to $0$, so that
$H$ is a self-adjoint operator and
$C^\infty_0(\R^N)$ is a core for $H$ (cf. \cite[Theorem~XIII.96]{rs4}). 
Furthermore, $H$ is bounded below.

\subsection{Approximating subspaces in  the quadratic projection 
method for the Schr\"odinger operator}

We have already established, in the abstract setting of Theorem~\ref{thm:lesh}
that, for any choice  of a subspace $\mathcal{L}\subset W^{2,2}(\R^N)$, the
eigenvalues of the matrix polynomial $P_\mathcal{L}(z)$ lying close  to the real
axis will be close to the spectrum of $H$ (and those ``far away'' from the real
axis don't matter). In  other words, the quadratic projection method does not
pollute. In order, however, to achieve a small error  and approximate as many
eigenvalues as possible, the choice of $\mathcal{L}$ (or of a sequence of such
spaces) is absolutely crucial, see Theorem~\ref{thm:approx}. Two main
difficulties here are the infinite geometry and the extra smoothness
requirements  needed for $Q_\mathcal{L}$  to make sense, see \eqref{eq:BAQ}.

Let $\Omega_{s}:=[-s,s]^N$. Let $\mathcal{M}_s:=W^{2,2}_0(\Omega_s)$ be a nested family of Sobolev spaces. 
Let $\mathcal{M}_{s,n}$, $n\in\N$, be a sequence of $n$-dimensional subspaces of $\mathcal{M}_s$. 
Let $\{\phi_{s,n,k}\}_{k=1}^n$ be a basis for $\mathcal{M}_{s,n}$. Set, for $j,k=1,\dots,n$,  
\begin{equation}\label{eq:matr_coeff}
\begin{aligned}
{}[B_{s,n}]_{j,k}&:=\int_{\Omega_s} \phi_{s,n,j}\phi_{s,n,k}\,,\\
[A_{s,n}]_{j,k}&:=\int_{\Omega_s} \nabla\phi_{s,n,j}\cdot\nabla\phi_{s,n,k}+V\phi_{s,n,j}\phi_{s,n,k}\,,\\
[Q_{s,n}]_{j,k}&:=\int_{\Omega_s} \Delta\phi_{s,n,j}\Delta\phi_{s,n,k}+2V\phi_{s,n,j}\Delta\phi_{s,n,k}
+V^2\phi_{s,n,j}\phi_{s,n,k}\,.
\end{aligned}
\end{equation}
and consider a quadratic $(n\times n)$-matrix polynomial
\begin{equation} \label{eq:P_sch}
P_{s,n}(z):=Q_{s,n}-2zA_{s,n}+B_{s,n}\,.
\end{equation}

Now, let $s_n$ be a monotone increasing unbounded sequence of positive real
numbers, let $\mathcal{L}_n = \mathcal{M}_{s_n,n}$, and let
$P_{n}(z)=P_{s_n,n}(z)$. Then Theorem~\ref{thm:approx} still holds  as long as
one can verify \eqref{eq:conv_cond} for $p=0,1,2$.

If the potential $V$ is sufficiently smooth, a natural choice of the basis
functions $\phi_{s,n,k}$ are  piecewise $C^2$ splines on $\Omega_s$ satisfying
$\phi|_{\partial\Omega_s}=\partial\phi/\partial n|_{\partial\Omega_s}=0$.
However, even for this simple choice, verifying  \eqref{eq:conv_cond} is still
highly technical, and we omit the details.

Even fixing \emph{both} parameters $s$ and $n$ and not imposing \emph{any}
condition  on $\mathcal{M}_{s,n}$ except $\mathcal{M}_{s,n}\subset
W^{2,2}_0(\Omega_s)$, still usually provides some useful information  about the
spectrum, with \emph{a posteriori} two-sided estimates: if
$\lambda_n\in\Spec(P_{s,n})$ and 
$\operatorname{dist}(\Re\lambda_n,\Spec_\ess(H))\ge|\Im\lambda_n|$, then there
exists $\lambda\in\Spec_\dsc(H)$ which lies in the the same spectral gap as
$\Re\lambda_n$. See Corollary \ref{cor:est2} for a sharper estimate.

On the other hand, to achieve approximation it is crucial that both parameters
$n$ and $s$ go to infinity in our choice of approximate spaces $\mathcal{L}_n$.
If we fix an arbitrarily large $s$ and let $n\to\infty$, then, though we still
do not have pollution  (unlike a standard projection method), neither
we have approximation.

\subsection{The quadratic matrix polynomial problem}

The quadratic projection method prescribes finding the spectrum $\Spec(P)$ of the
a quadratic matrix polynomial of the form
\[
   P(z)=Q+2z A+z^2 B,
\]
cf. Section~\ref{subsect:quad}. In applications, the matrix coefficients $Q$, $A$ and 
$B$ are expected to be sparse and real. They are always hermitean, so
$P(z)$ is a self-adjoint matrix polynomial in the sense of 
\cite{glr}.

The standard way of finding $\Spec(P)$ is to construct a suitable
companion linear pencil eigenvalue problem,
\begin{equation} \label{equ:pen}
    Lv = \mu K v 
\qquad \mathrm{for\ some\ } 0\not=v\in \mathcal{L}\oplus\mathcal{L}\,,
\end{equation}
such that $\mu\in \Spec(P)$ if and only if
\eqref{equ:pen} holds true. The coefficients, $L$, $K$,
of the companion form, $(L- z K)$, are
twice the size of the coefficients of $P(z)$. They are not unique. 
Two possible companion forms are given by:
\begin{equation*}
    L=\begin{pmatrix} 0 & N \\ -Q & A \end{pmatrix}
\qquad
    K=\begin{pmatrix} N & 0 \\ 0 & B \end{pmatrix}
\end{equation*}
and
\begin{equation*}
    L=\begin{pmatrix} -Q & 0 \\ 0 & N \end{pmatrix}
\qquad
    K=\begin{pmatrix} A & B \\ N & 0 \end{pmatrix}, 
\end{equation*}
where $N$ is a non-singular matrix. 

Different companion forms lead to different stability properties of the 
linear pencil problem to be solved once the matrices have been assembled. 
It is desirable finding a companion form that does not worsen the 
condition numbers of the original matrix polynomial spectral
problem. For a thorough account on this issue see \cite{hmt} and references 
therein.

\subsection{Examples}\label{subsect:ex}

\subsubsection*{One-dimensional example --- Gaussian perturbation of the Mathieu potential}

Let 
\begin{equation*}
N=1,\qquad\qquad V_p(x)=\cos(x), \qquad \qquad V_d(x)=-e^{-x^2},
\end{equation*}
and $H$ as in \eqref{eq:schr} with potential $V=V_p+V_d$. We now illustrate
how to implement the theoretical discussion carried out in the previous
sections to the study of $\Spec(H)$.

The essential spectrum of $H$ is determined by $V_p$.
It comprises an infinite number of non-intersecting bands $[\alpha_n,\beta_n]$
 of
absolutely continuous spectrum whose endpoints are determined by the Mathieu 
characteristic values \cite[\S 7.4]{in}. The approximate endpoints of the 
first five bands  
are given in Table~\ref{tab:tab1}.

\begin{table}[ht!]
\centerline{\begin{tabular}{|c|r|r|} 
\hline
$n$ & $\alpha_n$\hspace{.5cm}\ & $\beta_n$\hspace{.5cm}\ \\ 
\hline 
$1$ & $-0.378490 $ &  $-0.347670$ \\
$2$ & $0.594800$ & $0.918058$\\  
$3$ & $1.29317$ & $2.28516$\\
$4$ & $2.34258$ & $4.03192$\\
$5$ & $4.03530$ & $6.27082$\\
\hline
\end{tabular}}
\caption{Endpoints of the first five bands of the essential spectrum for
the Gaussian perturbation of the Mathieu potential. 
Cf. \cite{abst}.\label{tab:tab1}} 
\end{table}

Addition of the negative Gaussian potential yields a
non-empty discrete spectrum. By implementing the quadratic projection 
method \eqref{eq:P_sch}
into a finite element scheme, we detect
three eigenvalues of $H$ with high accuracy:
\[
  \lambda_1\approx -0\point 40961, \qquad
  \lambda_2\approx 0\point 37763, \qquad
  \lambda_3\approx 1\point 18216. \qquad
\]
The eigenvalue 
$\lambda_1$ is below the bottom of the essential spectrum, whereas
$\lambda_2$ and $\lambda_3$ lie in the first and the second gap, respectively.

\begin{figure}[ht!]
\centerline{\resizebox{12cm}{!}{\includegraphics*{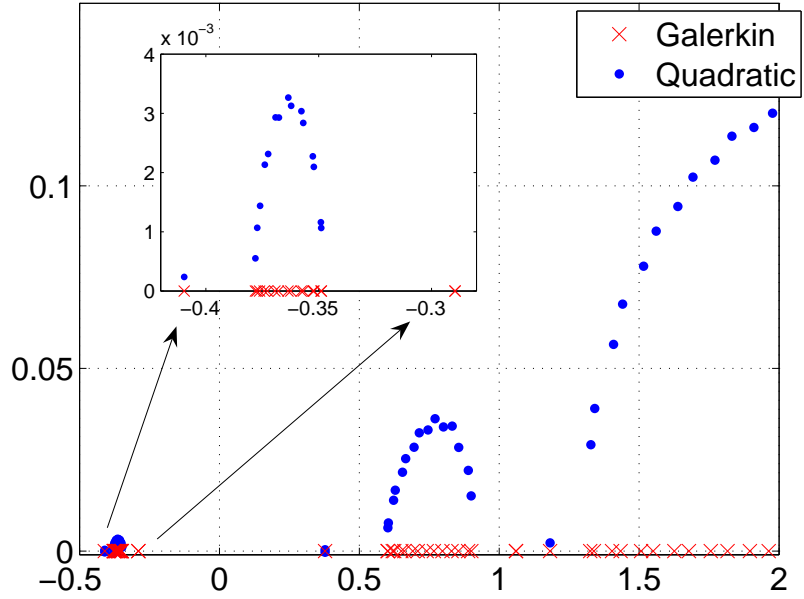}}}
\caption{The quadratic projection method vs the Galerkin projection method.  
Here $s=49$. Insert: zoom near $\lambda\approx-0.35$.\label{fig:fig1}}
\end{figure}

Figure~\ref{fig:fig1} illustrates the main ideas discussed in the previous sections. 
The spectrum of $P(z)$ is shown as blue
dots, while the eigenvalues of the standard Galerkin eigenvalue problem \eqref{eq:proj} are shown as red
crosses. The picture shows a narrow strip of the complex plane with the bottom
edge being the interval $[-0.5, 2]$. Note that there are eigenvalues
of $P(z)$ close to each of the eigenvalues $\lambda_1$,
$\lambda_2$ and $\lambda_3$. According to Corollary~\ref{cor:est1}, these eigenvalues are not spurious: the real
part of a complex number $z\in\Spec(P)$ is \emph{always} an
approximation of points in $\Spec (H)$, with a two-sided error estimate depending on $\Im z$. 
There are also eigenvalues of the linear problem \eqref{eq:proj} near to $\Spec_\dsc (H)$. These
eigenvalues also provide one-sided approximation (from above) for the
$\lambda_j$. However one should be careful when using the Galerkin
method, as spectral pollution may happen. For this particular set
of parameters there are two spurious eigenvalues: one near
$-0.3$ and the other near $1.1$.

\subsubsection*{Two-dimensional examples} 

We now consider a family of case studies with $N=2$. 
For $f\in W^{2,2}(\R^2)$, let
\begin{gather*}
   H_0 f(x,y)=-\Delta f(x,y)+(\cos(x)+\cos(y))f(x,y) \\
   H_1 f(x,y)=H_0 f(x,y)-c e^{-(x^2+y^2)}
   f(x,y) \\
H_2 f(x,y)=H_0 f(x,y)-c xe^{-(x^2+y^2)}
   f(x,y),  
\end{gather*}
where $c>0$.
A straightforward argument involving separation of variables shows that
\[
  \Spec (H_0) =\bigcup_{\lambda \in \Spec(\tilde{H})}
   \{\lambda+\mu\,:\, \mu\in \Spec(\tilde{H})\}  ,
\]
where $\tilde{H}=-\partial_x^2+\cos(x)$ is the one dimensional
Mathieu Hamiltonian. Furthermore, as both $H_1$ and $H_2$ are relatively compact 
perturbations of $H_0$, 
\[
   \Spec_\ess (H_1)=\Spec_\ess (H_2)=\Spec (H_0).
\]
An approximation of the endpoints of the bands
comprising the essential spectrum is given in Table~\ref{tab:tab2}.
Unlike the one-dimensional model, we now have a finite number of gaps.
Note that the perturbation associated to $H_1$ is radially
symmetric and sign definite, 
while the one associated to $H_2$ is sign indefinite and not radially symmetric.

\begin{table}[ht!]
\centerline{\begin{tabular}{|c|r|r|}
\hline
$n$ & $\alpha_n$\hspace{.5cm}\ & $\beta_n$\hspace{.5cm}\ \\
\hline 
$1$ & $-0.756978$ & $-0.695338$ \\
$2$ & $0.216310$ & $0.570389$ \\
$3$ & $0.914677$ & $\infty$ \\
\hline
\end{tabular}}
\caption{Endpoints of $\Spec_\ess(H_0)$.}
\label{tab:tab2}
\end{table}

With the quadratic projection method
we have been able to detect some discrete eigenvalues
of $H_1$ and $H_2$ for different values of the coupling constant $c$.
These results are presented in Table~\ref{tab:tab3}.
As we increase the value of $c$, eigenvalues of $H_1$ are 
moving from right to left. From the numerical results, the same seems to be true for eigenvalues of $H_2$.
Note that if an eigenvalue is close to an end-point of a band of the essential spectrum, the estimate
\eqref{eq:est1} does not allow us to distinguish between this eigenvalue and the end-point of the band --- 
thus the gaps in Table~\ref{tab:tab3}.

\begin{table}[ht!]
\centerline{\small \begin{tabular}{|l|r|r||l|r|r|}
\hline
\multicolumn{3}{|c||}{Eigenvalues of $H_1$}&\multicolumn{3}{c|}{Eigenvalues of $H_2$}\\
\hline
$c$ & $\lambda_1$\hspace{2cm}\ & $\lambda_2$\hspace{2cm}\ & $c$ & $\lambda_1$\hspace{2cm}\ & $\lambda_2$\hspace{2cm}\ \\
\hline 
 $5.0$ &  & $-0.09697\pm 3.39 \cdot 10^{-4}$   &$10$  & $0.1377\pm 6.61\cdot 10^{-3}$ & $0.7559\pm 7.07\cdot 10^{-2}$  \\
 $5.2$ &  & $-0.17133 \pm 2.03 \cdot 10^{-4}$  &  $11$  & $0.0865\pm 2.79\cdot 10^{-3}$ & $0.681\pm 1.05\cdot 10^{-1}$ \\
 $5.4$ &  & $-0.25255 \pm 1.45 \cdot 10^{-4}$  &  $12$  & $0.0115\pm 1.28\cdot 10^{-3}$ & \\
 $5.6$ &  & $-0.33905 \pm 1.44 \cdot 10^{-4}$  &  $13$  & $-0.09190\pm 8.19\cdot 10^{-4}$ &  \\
 $5.8$ &  & $-0.42902 \pm 1.71 \cdot 10^{-4}$  &  $14$  & $-0.22250\pm 7.77\cdot 10^{-4}$ &  \\
 $6.0$ & $-0.76946 \pm 2.85 \cdot 10^{-3}$      &   $-0.51978 \pm 2.52 \cdot 10^{-4}$  & $15$  & $-0.3730\pm 1.17\cdot 10^{-3}$ &   \\
 $6.2$ & $-0.78612\pm 8.43 \cdot 10^{-4}$ & $-0.60527 \pm 5.21 \cdot 10^{-4}$  & $16$  & $-0.5279\pm 1.97\cdot 10^{-3}$ &\\
\hline              
\end{tabular}}
\caption{Approximated eigenvalues of $H_1$ and $H_2$.} \label{tab:tab3}
\end{table}

Note that an eigenvalue $\lambda_1$ of the Hamiltonian $H_1$ is
below the
bottom of the essential spectrum for $c\gtrapprox 6$. 
The Galerkin method could actually be
implemented to approximate this eigenvalue. The quadratic projection, however,
works whether an eigenvalue is 
in a gap or not, and also provides a good approximation in this case.

\begin{figure}[ht!]
\centerline{\resizebox{9cm}{!}{\includegraphics*{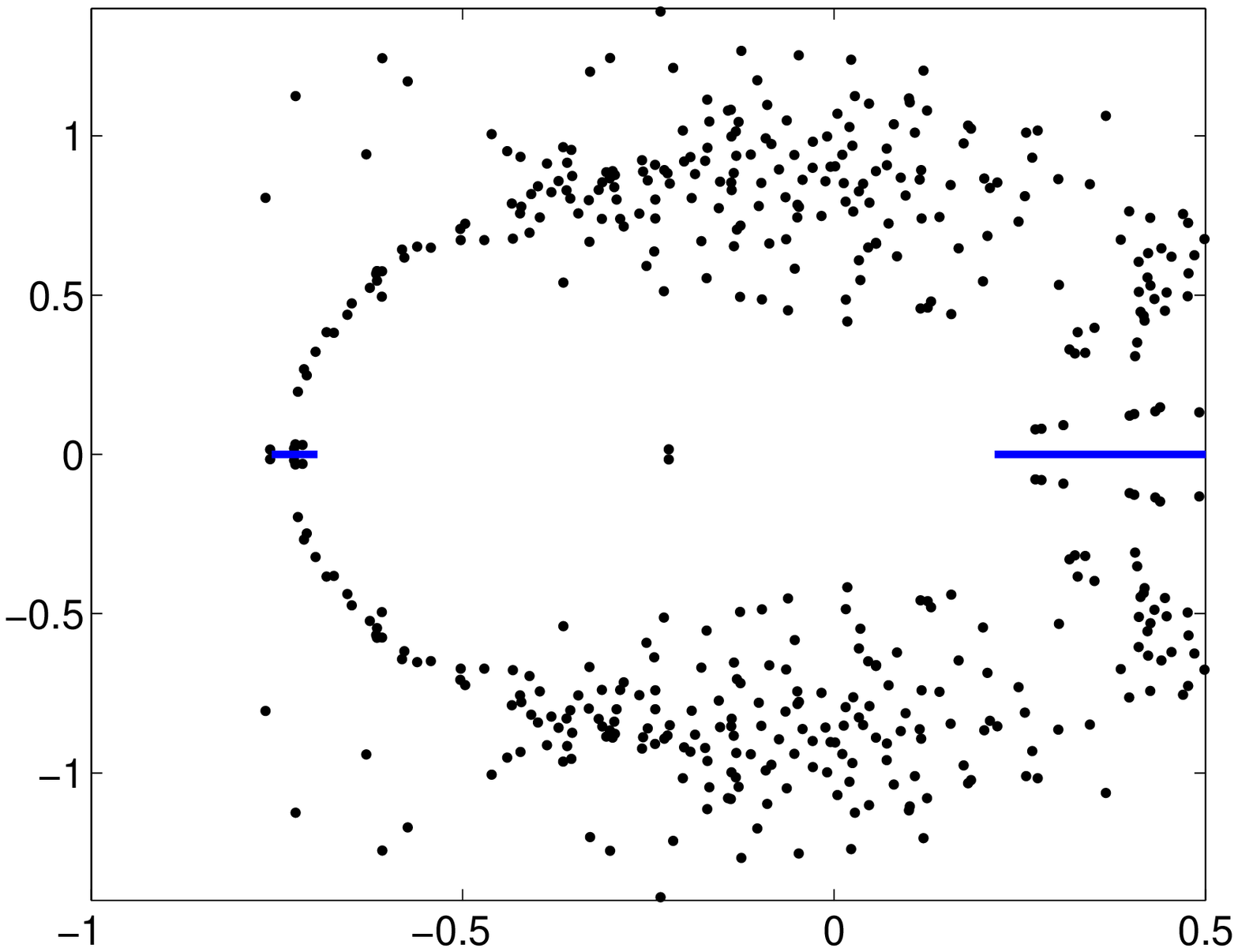}}
\hspace{-1cm}
\resizebox{9cm}{!}{\includegraphics*{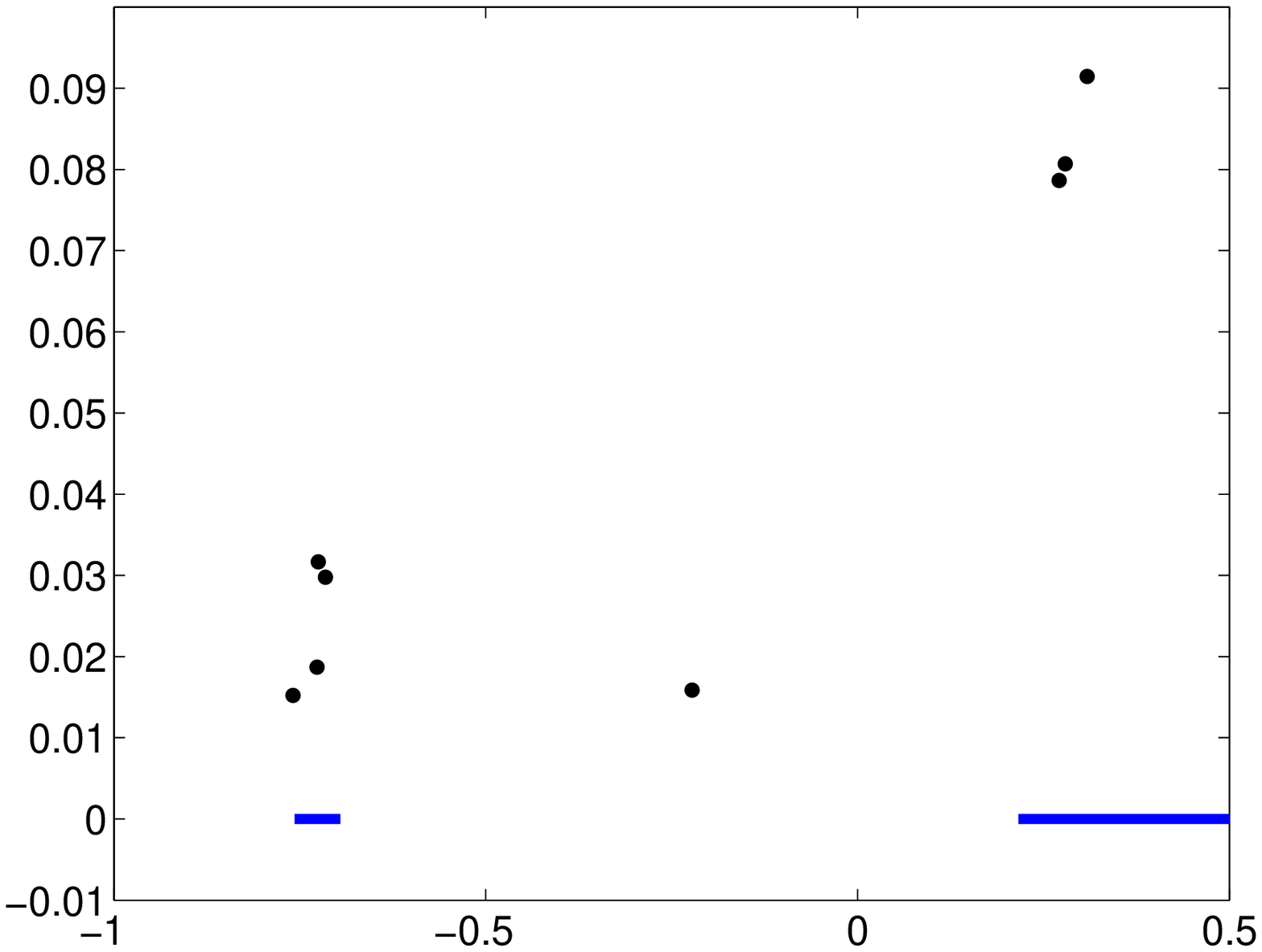}}}
\caption{The quadratic projection method for our two dimensional models.
Left: Typical output in the computation of $\Spec (P)$ for operator
$H_2$ when $c=14$ (here $s=60$). Right: zoom in the left picture
on a narrow strip near  the real line.  \label{fig:fig2}}
\end{figure}

In Figure~\ref{fig:fig2} we show the portion of
of $\Spec(P)$ lying in the box $[-1,1/2]\times[-3/2,3/2]$ 
for $H\equiv H_2$, $c=14$ and $s=60$. Corresponding
pictures for $H_1$ and other choices of $c$ and $s$ are qualitatively
similar. This graph clearly indicates approximation to an eigenvalue
$\lambda_1\approx -0.2225$ (see the right hand picture). A large portion
of $\Spec (P)$ forms an annular cloud around the spectral gap 
$(\beta_1,\alpha_2)$ and is sufficiently away from $\R$ to indicate
that there are no other eigenvalues in this gap. Note also that
some eigenvalue of $P(z)$ are close to $\Spec_\ess (H_2)$.

\section{Final remarks}

Other procedures exist for computing the eigenvalues of perturbed periodic 
partial differential operators such
as $H$, see \cite{do}. These include a method based on 
finding the eigenvalues of the matrix pencil problem
\begin{equation} \label{eq:proj_2}
   A_{s,n}u=\lambda  B_{s,n } u \qquad \mathrm{for\ some\ }
  u \in \mathcal{L}\setminus\{0\},
\end{equation}  
where the matrix coefficients are defined by \eqref{eq:matr_coeff}
(that is applying the projection method) for several values
of $s$ and $n$, and observing the dynamics of the eigenvalues
of \eqref{eq:proj_2} as $s$ increases. Some of the eigenvalues 
of \eqref{eq:proj_2} will be spurious and some will be close
to the true spectrum of $H$. The spurious eigenvalues will typically be
unstable as functions of the parameter $s$. 
The approximate eigenvalues close to the true spectrum of $H$ will be,
on the other hand, very stable. Thus, by increasing $s$, and tracking the
evolution of the eigenvalues of \eqref{eq:proj_2}, one would be able to
obtain some information about $\Spec(H)$.

This method, however, is quite inaccurate and it becomes useless when $N\geq 2$,
and we are interested in finding large eigenvalues. Furthermore,
it very much depends upon the choice of approximating 
subspaces $\mathcal{L}_s$. We are not aware of any
rigorous treatment of the effectiveness of this approach.

As the chosen subspaces $W_0^{2,2}(\Omega_s)$ are naturally nested for
increasing values of $s>0$ and they are all embedded in $W^{2,2}(\R^N)$,
every  point in $\Spec(H)$ is approximated (always from above) by the
spectrum of $(\overline{H\upharpoonright W_0^{2,2}(\Omega_s)})$.
Note that compactly supported functions form a core for 
the operator, and satisfy \emph{any} boundary condition
if the boundary is far enough away. Spectral pollution in the projection
method is a consequence of high eigenvalues of 
$(\overline{H\upharpoonright W_0^{2,2}(\Omega_s)})$ accumulating at
the bottom of the essential spectrum of $H$, and this effect is unavoidable.

We suggest using instead (or in addition to standard techniques), the quadratic projection
method, which never pollutes.

\subsection*{Acknowledgements} We are grateful to Marco Marletta for useful discussions and 
helpful advice.


\bibliographystyle{alpha}

\end{document}